\documentclass[a4paper,12pt]{article}
\usepackage{amsfonts,amssymb,epsfig}
\usepackage{amsmath}
\usepackage[latin1]{inputenc}
\newtheorem{thm}{Theorem}[section]
\newtheorem{cor}{Corollary}[section]
\newtheorem{lem}{Lemma}[section]
\newtheorem{pro}{Proposition}[section]
\newtheorem{dfn}{Definition}[section]
\newtheorem{rmk}{Remark}[section]
\newtheorem{expl}{Example}[section]
\newtheorem{nota}{Notation}[section]

\def\dessous#1\sous#2{\mathrel{\mathop{\kern0pt#2}\limits_{#1}}}

\newcommand{\al}{\alpha}
\newcommand{\be}{\beta}
\newcommand{\p}{\partial}
\newcommand{\n}{\nabla}
\newcommand{\pr}{\prime}
\newcommand{\mbb}{\mathbb}
\newcommand{\mf}{\mathfrak}
\newcommand{\mc}{\mathcal}
\newcommand{\nb}{\nonumber}
\newcommand{\ti}{\textit}
\newcommand{\beq}{\begin{eqnarray}}
\newcommand{\eeq}{\end{eqnarray}}
\newcommand{\bpro}{\begin{pro}}
\newcommand{\epro}{\end{pro}}
\newcommand{\blem}{\begin{lem}}
\newcommand{\elem}{\end{lem}}
\newcommand{\bdfn}{\begin{dfn}}
\newcommand{\edfn}{\end{dfn}}
\newcommand{\bcor}{\begin{cor}}
\newcommand{\ecor}{\end{cor}}
\newcommand{\bthm}{\begin{thm}}
\newcommand{\ethm}{\end{thm}}
\newcommand{\bex}{\begin{expl}}
\newcommand{\eex}{\end{expl}}
\newcommand{\brmk}{\begin{rmk}}
\newcommand{\ermk}{\end{rmk}}
\newcommand{\bnota}{\begin{nota}}
\newcommand{\enota}{\end{nota}}
\newcommand{\ben}{\begin{enumerate}}
\newcommand{\een}{\end{enumerate}}
\begin{document}
\title{$g$-Natural metrics on tangent bundles and Jacobi operators}
\author{S. Degla$^\ast$ and L. Todjihounde \thanks{Address: Institut de Math\'ematiques et de Sciences Physiques, 01 BP 613 Porto-Novo, Benin. E-mails: sdegla@imsp-uac.org ; leonardt@imsp-uac.org}}  
\date{}
\maketitle

\begin{abstract}

\noindent Let $(M,g)$ be a Riemannian manifold and $G$ a nondegenerate $g$-natural metric on its tangent bundle $TM$. In this paper we establish a relation between the Jacobi operators of $(M,g)$ and that of $(TM,G)$. \\ In the case of a Riemannian surface $(M,g)$, we compute explicitly the spectrum of some Jacobi operators of $(TM,G)$ and give necessary and sufficient conditions for $(TM,G)$ to be an Osserman manifold.
\\ \\
\noindent {\bf MSC 2000:} Primary 53B20,
53C07; Secondary 53A55, 53C25.
\\ \\
\noindent {\bf Key words}: $F$-tensor fields, $g$-natural metrics, Jacobi operators, Osserman manifolds.
\end{abstract}
\section*{0. Introduction}\label{s0}
 In \cite{AS1} the authors introduced $g$-natural metrics on the tangent bundle $TM$ of a Riemannian manifold $(M,g)$ as metrics on $TM$ which come from $g$ through first order natural operators defined between the natural bundle of Riemannian metrics on $M$ and   the natural bundle of $(0,2)$-tensors fields on the tangent bundles. Classical well-known metrics like Sasaki metric (cf. \cite{Sa} , \cite{Do}) or Cheeger-Gromoll metrics (cf. \cite{CG} , \cite{MT}) are examples of natural metrics on the tangent bundle. 
 By associating the notion of $F$-tensors fields they got a characterization of $g$-natural metrics  on $TM$ in terms of the basis metric $g$ and of some functions defined on the set of positive real numbers, and obtained  necessary and sufficient conditions for $g$-natural metrics to be either nondegenerate or Riemannian (see \cite{KMS} for more details on natural operators and $F$-tensors fields).
\\
\noindent Some geometrical properties of $g$-natural metrics are inherited from the basis metric $g$ and conversely (cf. \cite{AS1}, \cite{AS2}, \cite{DET}, \cite{KS2}). We will investigate in this paper the property of being Osserman which is closely related to the  spectrum of Jacobi operators. \\Recall that for a tangent vector $X\in T_xM$ with $x\in M$, the  Jacobi operator $J_X$ is defined as the linear self-adjoint map \\
$\begin{array}{llcl}     &       &           &    \\  
                    J_X: & T_xM &\rightarrow &T_xM\\
                          &  Y   & \mapsto & J_X(Y):=R(X,Y)X \,,
\end{array}$ where $R$ denotes the Riemannian curvature operator of $(M,g)$. Osserman manifolds are defined as follows:

\bdfn
1. Let $x\in M$. $(M,g)$ is  Osserman  at $x$ if, for any unit tangent vector $X \in T_xM$, the eigenvalues of the Jacobi operator $J_X$ do not depend on $X$. \\ \\
2. $(M,g)$ is pointwise Osserman  if it is Osserman 
at any point of $M$. \\ \\
3. $(M,g)$ is globally Osserman manifold if, for any point $x \in M$ and any unit tangent vector $X \in T_xM$,  
 the eigenvalues of the Jacobi operator $J_X$ depend neither on $X$ nor on $x$.  
\edfn
Globally Osserman manifolds are obviously pointwise Osserman.
\brmk  \label{l1s3}
For any point $x \in M$ the map defined on $T_xM$ by $X \longmapsto J_X$ satisfies the identity $J_{\lambda X} = \lambda^2 J_X\;,\;\forall \;\lambda \in {\mathbb R}$. So the spectrum of $J_{\lambda X}$ is, up to the factor $\frac{1}{\lambda^2}$, the same that of $J_X$.
Thus $(M,g)$ is Osserman at $x \in M$ if and only if for any vector $X\in T_xM$ with $X\neq 0$ and for any eigenvalue $\lambda(X)$ of $J_X$, the quotient  $\frac{\lambda(X)}{g(X,X)}$ does not depend on $X$.
\ermk 

\noindent Flat manifolds or locally symmetric spaces of rank one are examples of global\-ly Osserman manifolds since the local isometry group acts transitively on the unit tangent bundle, and hence the eigenvalues of the Jacobi operators are constant on the unit tangent bundle.\\
 Osserman conjectured that the converse holds; that is all Osserman manifolds are locally symmetric of rank one. The Osserman conjecture has been proved in many special cases (cf. \cite{Ch}, \cite{Ni}, \cite{Os}, \cite{SV}).
\\
Using the fact that $(M,g)$ is totally geodesic in $(TM,G)$ (cf. \cite{AS1}) we show that any eigenvalue of a Jacobi operator of $(M,g)$ is an eigenvalue  of some Jacobi operator of its $g$-natural  tangent bundle $(TM,G)$. Furthermore, we investigate the
Jacobi operators  of $g$-natural metrics on tangent bundles of Riemannian surfaces, and we compute explicitly their spectrums. Then 
we establish  necessary and sufficient conditions for 
$g$-natural tangent bundles of Riemannian surfaces to be Osserman manifolds. 
\section{Preliminaries}\label{s1}
Let $(M,g)$ be a Riemannian manifold and $\n$ the Levi-Civita 
connection of $g$. 
The tangent space of $TM$ at a point
$(x,u)\in TM$ splits into the horizontal and vertical subspaces
with respect to $\n$ :
$$ 
T_{(x,u)}TM = H_{(x,u)}M\oplus V_{(x,u)}M\; .
$$

\noindent A system of local coordinates $\left(U\, ;\, x_i,\; i=1,\cdots,m\right)$ in $M$ induces on $TM$ a system of local coordinates $\left(\pi^{-1}(U)\, ;\, x_i,u^i,\; i=1,\cdots,m\right)$.\\ Let 
$X=\sum_{i=1}^m X^i \frac{\p}{\p x_i}$ be the local expression in $U$ of a vector field $X$ on $M$. Then, the horizontal lift $X^h$ and the vertical lift $X^v$ of $X$ are given, \\ with respect to the induced coordinates, by :
\beq
X^h & = &\sum_iX^i\frac{\partial}{\partial x_i}-
\sum_{i,j,k}\Gamma_{jk}^iu^jX^k\frac{\partial}{\partial u^i}\quad\mbox{ and }\\
X^v & = & \sum_{i}X^i\frac{\partial}{\partial u^i},
\eeq
where the $\Gamma_{jk}^i$ are the   Christoffel's symbols defined by  $g$.
\\
\noindent Next, we  introduce some notations which will be used to describe\\ vectors
obtained from lifted vectors by basic operations on $TM$. Let $T$ be \\ a tensor
field of type $(1,s)$ on $M$. If $X_1,X_2,\cdots,X_{s-1}\in T_xM,$ then \\
$h\{T(X_1,\cdots, u,\cdots,X_{s-1})\}$ and $v\{T(X_1,\cdots, u,
\cdots,X_{s-1})\}$) are horizontal and vertical vectors repectively  at the point $(x,u)$ which are defined by:
$$ h\{T(X_1,\cdots, u,\cdots,X_{s-1})\}=\sum u^\lambda
\left(T(X_1,\cdots, \frac{\p}{\p x_\lambda}_{|x}\,,\cdots,X_{s-1})\right)^h $$
$$
v\{T(X_1,\cdots, u,\cdots,X_{s-1})\}=\sum u^\lambda
\left(T(X_1,\cdots, \frac{\p}{\p x_\lambda}_{|x}\,,\cdots,X_{s-1})\right)^v \;.$$
In particular, if $T$ is the identity tensor of type $(1,1)$, then we
 obtain the geodesic flow vector field at $(x,u)$,
 $\xi_{(x,u)} =\sum_{\lambda}u^{\lambda}\left(\frac{\p}{\p x_\lambda}\right)_{(x,u)}^h$, and the cano\-nical vertical vector at $(x,u)$,
 $\mc{U}_{(x,u)}=
 \sum_{\lambda}u^\lambda\left(\frac{\p}{\p x_\lambda}\right)_{(x,u)}^v$.
\\ 
Also  $h\{T(X_1,\cdots, u,\cdots,u,\cdots,X_{s-t})\}$
and \\
$v\{T(X_1,\cdots, u,\cdots,u,\cdots,X_{s-t})\}$ are defined by 
similar way. 
\\
Let us introduce the notations
\beq
h\{T(X_1,\cdots,X_s)\} =: T(X_1,\cdots,X_s)^h
\eeq
 and 
\beq
 v\{T(X_1,\cdots, X_s)\}=: T(X_1,\cdots,X_s)^v\; .
\eeq
Thus $h\{X\}=X^h$ and $v\{X\}=X^v$, for each vector field $X$ on $M$.
\\
From the preceding quantities, one can define vector fields on $TU$ in the following way: If $u=\sum_iu^i\left(\frac{\p}{\p x_i}\right)_x$ is a given point in $TU$ and $X_1,\cdots,X_{s-1}$ are vector fields
on $U$, then we denote by 
$$
h\{T(X_1,\cdots, u,\cdots,X_{s-1})\}\quad
(\mbox{respectively}\quad
v\{T(X_1,\cdots, u,\cdots,X_{s-1})\}) $$
the horizontal (respectively vertical) vector field on $TU$ 
defined by
$$
h\{T(X_1,\cdots, u,\cdots,X_{s-1})\}=\sum_\lambda u^\lambda
T(X_1,\cdots,\frac{\p}{\p x_\lambda}\,,\cdots,X_{s-1})^h
$$
$$
(\mbox{ resp.}\quad
v\{T(X_1,\cdots, u,\cdots,X_{s-1})\}=\sum_\lambda u^\lambda
T(X_1,\cdots,\frac{\p}{\p x_\lambda}\,,\cdots,X_{s-1})^v
\, ) .
$$
\noindent Moreover, for vector fields $X_1,\cdots,X_{s-t}$ on $U$, where $s\, ,\, t\in \mbb{N}^\ast\, (s>t) $,  the\\
vector fields $h\{T(X_1,\cdots, u,\cdots, u, \cdots,X_{s-t})\}$ and \\
$v\{T(X_1,\cdots, u,\cdots,u, \cdots,X_{s-t})\}$ on $TU$, are defined by similar way.
\\
Now, for $(r,s)\in \mbb{N}^2$, we denote by 
$\pi_M:\ TM\rightarrow M$ the natural  projection and $F$ the natural bundle 
defined by  
\beq
FM&=&\pi_M^*(\underbrace{T^*\otimes \cdots\otimes T^*}_{\mbox{$r$ times}}\otimes
\underbrace{T\otimes \cdots\otimes T}_{\mbox{$s$ times}})M\rightarrow M, \\ \nb
Ff(X_x,S_x)&=&(Tf.X_x,(T^*\otimes \cdots\otimes T^*\otimes 
T\otimes \cdots\otimes T)f.S_x) 
\eeq
for $x \in M$, $X_x \in T_xM$, $S \in (T^*\otimes \cdots\otimes T^*\otimes T\otimes \cdots\otimes T)M$ and any local diffeomorphism $f$ of $M$. 
\\
We call the sections of the canonical projection
$FM\rightarrow M$ $F$-tensor fields of type $(r,s)$.
 So, if we denotes the product of fibered  manifolds by $\oplus$, then the  $F$-tensor fields are mappings\\  
$A:\ TM\oplus \underbrace{TM\oplus\cdots \oplus TM}_{\mbox{$s$ times}} \rightarrow \underset{x\in M}{\sqcup}\otimes^rT_xM$ which are   linear in the last\\
$s$ summands and such that  $\pi_2\circ A=\pi_1$, where $\pi_1$ and
$\pi_2$ are respectively  the  natural projections of the source and target fiber bundles of $A$.
For $r=0$ and $s=2$, we obtain the classical notion of $F$-metrics. So,
$F$-metrics are mappings $TM\oplus TM\oplus TM\rightarrow \mbb{R}$ 
which are linear in the second and  the third arguments.
\\
Moreover let us fix $(x,u)\in TM$ and a system of normal coordinates \\
$S:= (U\, ;\, x_i\, ,i=1,\cdots, m)$ of $(M,g)$ centred at $x$. Then we can define on $U$ the vector field $\mathbf{U} := \sum_iu^i\frac{\p}{\p x_i}$,
where $(u^1,\cdots,u^m)$ are the coordinates of $u\in T_xM$ with respect to its basis $ (\frac{\p}{\p x_i}_{|x}; \, i=1,\cdots,m )$.
\\
Let $P$ be an $F$-tensor field of type $(r,s)$ on $M$. Then on $U$, we can define an $(r,s)$-tensors field $P_u^S$ (or $P_u$ if there is no risk of confusion) associated to $u$ and $S$ by
\beq
P_u(X_1,\cdots,X_s) := P(\mathbf{U}_z; X_1,\cdots,X_s)\, ,
\eeq
for all $(X_1,\cdots,X_s)\in T_zM,\; \forall z\in U$.
\\
On the other hand, if we fix $x\in M$ and $s$ vectors $X_1,\cdots,X_s$ in 
$T_xM$, then we can define a $C^\infty$-mapping
$P_{(X_1,\cdots,X_s)}: T_xM\rightarrow \otimes^rT_xM$, associated to 
$(X_1,\cdots,X_s)$ by 
\beq
P_{(X_1,\cdots,X_s)}(u) :=P(u;\, X_1,\cdots,X_s)\, ,
\eeq
for all $u\in T_xM$.
\\
Let $s>t$ be two non-negative integers, $T$ be a $(1,s)$-tensor field on $M$ and $P^T$ be an $F$-tensor field of type $(1,t)$ of the form
\beq
P^T(u; X_1,\cdots, X_t) = T(X_1,\cdots, u,\cdots, u,\cdots, X_t) ,
\eeq
for all $(u; X_1,\cdots, X_t)\in TM\oplus\cdots\oplus TM $, i.e., $u$ appears
$s-t$ times at positions $i_1,\cdots,i_{s-t}$ in the expression of $T$. Then
\begin{itemize}
\item[{\bf -}] $P_u^T$ is a $(1,t)$-tensor field on a neighborhood $U$ of
$x$ in $M$,\\ for all $u\in T_xM$~;
\item[{\bf -}] $P_{(X_1,\cdots,X_t)}^T$ is a $C^\infty$-mapping 
$T_xM\rightarrow T_xM$, for all $X_1,\cdots,X_t$ in $T_xM$. 
\end{itemize}
Furthermore, it holds 
\blem\cite{AS2}\label{l1s3}
\begin{enumerate}
\item[1)] The covariant derivative of $P_u^T$, with respect to the Levi-Civita connection of $(M,g)$ is given by :
\beq
\left(\n_X P_u^T\right)(X_1,\cdots,X_t) = (\n_X T)
(X_1,\cdots,u,\cdots,u, X_t),
\eeq
for all vectors $X,X_1,\cdots,X_t$ in $T_xM$, where $u$ appears at
positions \\ $i_1,\cdots,i_{s-t}$ in the right-hand side of the 
preceding formula.
\item[2)] The differential of $P_{(X_1,\cdots,X_t)}^T$ at $u\in T_xM$, is
given by :
\beq
d\left(P_{(X_1,\cdots,X_t)}^T\right)_u(X) &=&
T(X_1,\cdots,X,\cdots, u,\cdots,X_t) +\cdots\\ \nb 
 & &   + T(X_1,\cdots,u,\cdots, X,\cdots,X_t),
\eeq
for all $X\in T_xM$.
\end{enumerate}
\elem
\section{$g$-natural metrics on tangent bundles}\label{s2}

\bdfn
Let $(M,g)$ be a Riemannian manifold. A $g$-natural metric on the tangent bundle of $M$ is a metric on $TM$ which is the image of $g$ by a first order natural operator defined from  the natural bundle of Riemannian metrics $S_+^2T^*$ on $M$ into the natural bundle of $(0,2)$-tensor fields $(S^2T^*)T$ on the tangent bundles (cf. \cite{AS1} , \cite{AS2}).\\ \\
Tangent bundles equipped with $g$-natural metrics are called $g$-natural tangent bundles. 
\edfn

\noindent The following result gives the classical expression of $g$-natural metrics:
\bpro\cite{AS1}\label{p1s2}
Let $(M,g)$ be  a Riemannian manifold and $G$  a \\  $g$-natural metric on $TM$. There exists six smooth functions  \\ $\al_i,\
 \be_i:\mathbb{R}^+\to \mathbb{R},\ i=1,2,3,$ such that for any $x\in M$ and all  vectors $u,\ X,\ Y\in T_xM$, we have
\beq
\left\{
\begin{array}{lcl}
G_{(x,u)}\left(X^h,Y^h\right)
   &=& (\al_1+\al_3)(t)g_x(X,Y)\\ \nb
  &&+(\be_1+\be_3)(t)g_x(X,u)g_x(Y,u),\\ \nb
   & & \\ \nb
G_{(x,u)}\left(X^h,Y^v\right)
   &=& \al_2(t)g_x(X,Y)
      +\be_2(t)g_x(X,u)g_x(Y,u),\\ \nb
       & & \\ \nb
G_{(x,u)}\left(X^v,Y^h\right)
   &=& \al_2(t)g_x(X,Y)
         +\be_2(t)g_x(X,u)g_x(Y,u),\\ \nb
  & & \\ \nb
G_{(x,u)}\left(X^v,Y^v\right)
   &=& \al_1(t)g_x(X,Y)
         +\be_1(t)g_x(X,u)g_x(Y,u), \nb
\end{array}
\right.
\eeq
where $t=g_x(u,u)$, $X^h$ and $X^v$ are  respectively the  horizontal lift and the vertical lift  of the  vector $X\in T_xM$ at the  point  $(x,u)\in TM$.
\epro

\bnota\label{n1s2}\
\begin{itemize}
\item $\phi_i(t)=\al_i(t)+t\be_i(t)$, $i=1,2, 3$,
\item $\al(t)=\al_1(t)(\al_1+\al_3)(t)-\al_2^2(t)$,
\item $\phi(t)=\phi_1(t)(\phi_1+\phi_3)(t)-\phi_2^2(t)$,
\end{itemize}
for all $t\in \mbb{R}^+$.
\enota
For a $g$-natural metric to be nondegenerate or Riemannian, there are some conditions to be satisfied by the functions $\al_i$ and $\be_i$ of Proposition 2.1. It holds:
\bpro\cite{AS1}\label{p2s2}
A $g$-natural metric $G$ on the tangent bundle of a\\ Riemannian  manifold $(M,\ g)$ is :
\begin{enumerate}
\item[(i)]  nondegenerate if and only if the functions $\al_i,\, \be_i,\,
i=1,2,3$ defining $G$ are such that
\beq
\al(t)\phi(t)\neq 0
\eeq
for all $t\in \mbb{R}^+$.
\item[(ii)]  Riemannian if and only if the functions  
$\al_i,\, \be_i,\,i=1,2,3$ defining $G$,  satisfy the  inequalities
\beq
\left\{
\begin{array}{ll}
\al_1(t)>0,& \phi_1(t)>0,\\
\al(t)> 0, &\phi(t)>0,
\end{array}
\right.
\eeq
for all $t\in \mbb{R}^+$.

For $dim\,M=1$, this system  reduces to $\al_1(t)>0$
and $\al(t)>0$, for all $t\in \mbb{R}^+$.
\end{enumerate}
\epro

\noindent Before giving the formulas relating both Levi-Civita connexions $\n$ of $(M,g)$ and $\bar{\n}$ of $(TM , G)$ let us introduce the following notations: 

\bnota\label{n2s2}
For a Riemannian manifold $(M,g)$,  we set :
\beq
\begin{array}{lcl}
T^1(u;X_x,Y_x)=R(X_x,u)Y_x, & & T^2(u;X_x,Y)=R(Y_x,u)X_x\, ,\\
 T^3(u;X_x,Y_x)=g(R(X_x,u)Y_x,u)u , & & T^4(u;X_x,Y_x)=g(X_x,u)Y_x\\
 T^5(u;X_x,Y_x)=g(Y_x,u)X_x, & & T^6(u;X_x,Y_x)=g(X_x,Y_x)u,    \\
 T^7(u;X_x,Y_x)=g(X_x,u)g(Y_x,u)u. & & 
\end{array}
\eeq
where $(x,u)\in TM\, $,  $\, X_x,Y_x\in T_xM\, $ and $\, R\, $ is the Riemannian
curvature of $g$.
\enota
The $g$-natural metric $G$ being defined by the functions $\al_i ,  \be_i$ of Proposition 2.1, it holds:

\bpro\cite{DET}\label{p1s2}
Let $(x,u)\in TM$ and $X,Y\in \mf{X}(M)$, we have
\beq\label{e1s2}
\left(\bar{\n}_{X^h}Y^h\right)_{(x,u)} &=&\left(\n_XY\right)_{(x,u)}^h
              + h\{A(u;X_x,Y_x)\}+v\{B(u;X_x,Y_x)\}\quad   \label{h1}\\ 
\left(\bar{\n}_{X^h}Y^v\right)_{(x,u)} &=&\left(\n_XY\right)_{(x,u)}^v 
              + h\{C(u;X_x,Y_x)\}+v\{D(u;X_x,Y_x)\}\quad \label{h2} \\ 
\left(\bar{\n}_{X^v}Y^h\right)_{(x,u)} &=& h\{C(u;Y_x,X_x)\}+v\{D(u;Y_x,X_x)\}
\label{h3}\\
\left(\bar{\n}_{X^v}Y^v\right)_{(x,u)} &=& h\{E(u;Y_x,X_x)\}+v\{F(u;Y_x,X_x)\}
\label{h4}
\eeq
where
$P(u;X_x,Y_x)=\sum_{i=1}^8f^P_i(|u|^2)T^i(u;X_x,Y_x)\;,\; \mbox{for }
P=A,B,C,D,E,F$, and the functions $f_i^P$ defined as in \cite{DET}.
\epro
\noindent It has been notified by the authors in \cite{AS1} that,  
 the Riemannian manifod $(M,g)$, considered as an embedded submanifold in its $g$-natural tangent bundle $(TM , G)$ by the null section, is always totally geodesic. \\ Indeed the null section $S_0$ of $\mf{X}(M)$  is defined by
\beq
\begin{array}{lccl}
S_0:& M &\rightarrow & TM\\
    & x & \mapsto & (x,\ 0_x)\;,
\end{array}
\eeq
which determines an embedding of $M$ in $TM$. \\ Its differential at any point $x\in M$ is given by 
\beq
\begin{array}{lccl}\label{e2s2}
d{S_0}_{|x}:& T_xM & \rightarrow& T_{(x,\, 0_x)}TM\\
    & X_x& \mapsto &  X_{(x,0_x)}^h
    \end{array}\;.
    \eeq
Then according to (\ref{e1s2}) and (\ref{e2s2}) we have
\beq
\bar{\n}_{S_*X}S_*Y = \bar{\n}_{X^h \circ S}(Y^h \circ S) = S_* (\n_XY)\;,
\eeq
 for all $X,Y \in \mf{X}(M)$. \\ \\
Thus from the relation $(20)$ it holds: 
\bpro\cite{AS1}\label{p2s2}
Any Riemannian manifold $(M,g)$ is totally geodesic in its tangent bundle $TM$ equipped with a non-degenerate $g$-natural metrics $G$.
\epro
\brmk
 If  $G$ is nondegerate then the orthogonal of $S_0(M) \equiv M$ 
in $(TM, G)$ is given by
\beq
T_xM^{\perp_G} = \{H_{(x,0_x)}^h + V_{(x,0_x)}^v\in T_{(x,\, 0_x)}TM / \\ \nb 
H,V\in T_xM \mbox{ and }(\al_1+\al_3) H + \al_2 V = 0_x  \} ,
\eeq
where the functions $\al_i,\; i=1,2,3$ are evaluated at  $0$.
\ermk
\section{Jacobi operators and Osserman $g$-natural tangent bundles}
\label{s3}
\noindent In the above section, we mentioned that $(M,g)$ is totally geodesic in $(TM,G)$. By using this observation we get the following result:
 
\bpro
Assume that $\dim M\geq 2$, and $x\in M$.  If $\lambda$ is an eigenvalue of a Jacobi operator $J_X$ for $X\in S(T_xM)$ then $\lambda$ is an eigenvalue of the Jacobi operator $\bar{J}_{X_{(x,0_x)}^h}$ of $G$ at the point $(x,0_x)\in TM$.
\epro

\noindent \ti{Proof}\\
Let us choose in $(T_xM, g_x)$ an orthonormal basis $(X_1,\cdots X_m)$ such as \\$X_1=X$, and an orthonormal basis $(V_1,\cdots V_m)$ in $T_xM^{\perp_G}$. \\ Then $(X_1^h|_{(x,0_x)},\cdots X_m^h|{(x,0_x)},V_1,\cdots, V_m)$ is an orthogonal basis of $T_{(x,0_x)}TM$. Since
$(M,g)$ is totally geodesic in $(TM,G)$ and $\bar{J}_{X^h_{(x,0_x)}}$ is self-adjoint, the matrix of
$\bar{J}_{X_{(x,0_x)}^h}$ in this basis has the form
$\left(\begin{array}{ccc}J_1 & & 0\\ 0 & & J_2 \end{array}\right)\;,$ where $J_1$ is the matrix of
$J_X$ in the basis  $(X_1,\cdots X_m)$,  and $J_2$ is a square matrice of order $m$. Thus if $\lambda$ is an eigenvalue of $J_X$ then $\lambda$ is an eigenvalue of $\bar{J}_{X_{(x,0_x)}^h}$.
$\hfill{\square}$
\\ \\
As a collorary we have:

\bcor
 If $(TM,G)$ is pointwise Osserman manifold (respectively \\ globally Osserman manifold) then the same holds for $(M,g)$. 
\ecor
\noindent We shall now give the explicit expression of $\bar{J}$ in termes of the Levi-Civita connexion $\n$ and the curvature tensor $R$ of $(M,g)$ and some $F$-tensors on $M$. 
\\ \\
Let $(x,u)\in TM$ and $\bar{X} = H_{(x,u)}^h+V_{(x,u)}^v \in T_{(x,u)}TM$, with
$H\in T_xM$ and $V\in T_xM$. In the following we give an expression of the Jacobi operator $\bar{J}_{\bar{X}}$ of $(TM,G)$.
Firstly, let us consider the following $F$-tensors which are defined in terms of the
$F$-tensors $A,\, B,\, C,\, D,\, E,\mbox{ and } F$ of Proposition \ref{p1s2} such
as we have at any point $x\in M$:
\beq
P_{(A,B,C)}^1(u; H,Y,V)&=& B(u;H,A(u;Y,H))-B(u;Y,A(u;H,H))\\ \nb
                     &&+B(u;H,C(u;Y,V))- 2B(u; Y,C(u;H,V))
\eeq
\beq
P_{(A,C,D,F)}^2(u; H,Y,V)&=& D(u;A(u;Y,H),V) + D(u;C(u;Y,V),V)\, \\ \nb
                         & & -D(u;Y,F(u;V,V)),
\eeq
\beq
P_{(B,C,D)}^3(u; H,Y,V)&=& C(u;H,B(u;Y,H))-C(u;Y,B(u;H,H))\\ \nb
                     &&+C(u;H,D(u;Y,V))- 2C(u; Y,D(u;H,V)),
\eeq
\beq
P_{(A,B,D,E,F)}^4(u; H,Y,V)&=& F(u;V,B(u;Y,H))\\ \nb
                           &&+F(u;V,D(u;Y,V))-A(u;Y,E(u;V,V)),
\eeq
\beq
P_{(C,E)}^5(u; H,Y,V)= C(u;H,R(H,Y)u)+  E(u;R(H,Y)u,V),
\eeq
\beq
P_{(A,C)}^6(u; H,Y,V)= d\left(A_{(Y,H)}\right)_u(V)+ d\left(C_{(Y,V)}\right)_u(V) ,
\eeq
\beq
Q_{(A,B,C, D,F)}^1(u; H,Y,V)&=& A(u;H,C(u;H,Y))+A(u;H,E(u;Y,V))\\ \nb
                     &&-F(u;Y,B(u;H,H))- 2F(u; Y,D(u;H,V)),
\eeq
\beq
Q_{(D,E,F)}^2(u; H,Y,V) &=& E(u;V,D(u;H,Y))
                           +E(u;V,F(u;Y,V))\\ \nb
                     &&-E(u;Y,F(u;V,V)),
\eeq
\beq
Q_{(A,B,C)}^3(u; H,Y,V)&=&D(u;C(u;H,Y),V)-2D(u;C(u;H,V),Y)\\ \nb
                     &&+D(u;E(u;Y,V),V)- D(u; E(u;V,V),Y) ,
\eeq
\beq
Q_{(A,C,D,F)}^4(u; H,Y,V) &=& C(u;H,D(u;H,Y))\\ \nb
                           &&+C(u;H,F(u;Y,V))
                     -C(u;A(u;H,H),Y) ,
\eeq
\beq
Q_{(C)}^5(u; H,Y,V)= d\left(C_{(H,Y)}\right)_u(V)-2 d\left(C_{(H,V)}\right)_u(Y),
\eeq
\beq
Q_{(A,E)}^6(u; H,Y,V)&=& d\left(E_{(Y,V)}\right)_u(V)
- d\left(E_{(V,V)}\right)_u(Y)\\ \nb
 && -d\left(A_{(H,H)}\right)_u(Y),
\eeq
for all $u,H,Y,V\in T_xM$.
\\ \\
\noindent The Jacobi operator $\bar{J}_{\bar{X}}$ is then determined by
\beq\label{e1s3}
\bar{J}_{\bar{X}}\left(Y^h\right)&=& h\{ R(H,Y)H + 
           \left[\left(\n_HA_u\right)(Y,H)-\left(\n_YA_u\right)(H,H)\right]\\ \nb
 && +\left[\left(\n_HC_u\right)(Y,V)-\left(\n_YC_u\right)(H,V)\right]\\ \nb
 && -\left(\n_YC_u\right)(H,V)-\left(\n_YE_u\right)(V,V)\\ \nb
  &&+ 
P_{(A,A,C)}^1(u; H,Y,V)+
P_{(A,C,C,F)}^2(u;H,Y,V)+
P_{(B,C,D)}^3(u; H,Y,V)\\ \nb
&&+P_{(A,B,D,E,E)}^4(u; H,Y,V)+
P_{(C,E)}^5(u; H,Y,V)+
P_{(A,C)}^6(u; H,Y,V)
\} \\ \nb
 & +& \\
&& v\{R(H,Y)V+ \left[\left(\n_HB_u\right)(Y,H)-\left(\n_YB_u\right)(H,H)\right]\\ \nb
 && +\left[\left(\n_HD_u\right)(Y,V)-\left(\n_YD_u\right)(H,V)\right]\\ \nb
 && -\left(\n_YD_u\right)(H,V)-\left(\n_YF_u\right)(V,V)\\ \nb
  &&+ 
P_{(A,B,C)}^1(u; H,Y,V)+
P_{(A,C,D,F)}^2(u;H,Y,V)+
P_{(B,D,D)}^3(u; H,Y,V)\\ \nb
&&+P_{(B,B,D,E,F)}^4(u; H,Y,V)+
P_{(D,F)}^5(u; H,Y,V)+
P_{(B,D)}^6(u; H,Y,V)
\}
\eeq
and 
\beq\label{e2s3}
\bar{J}_{\bar{X}}\left(Y^v\right)&=&h\{ \left(\n_HC_u\right)(H,Y)
+\left(\n_HE_u\right)(Y,V)\\ \nb
&& +Q_{(A,B,C, D,E)}^1(u; H,Y,V)+
Q_{(D,E,F)}^2(u; H,Y,V)+
Q_{(C,C,E)}^3(u; H,Y,V)\\ \nb
&&+Q_{(A,C,D,F)}^4(u; H,Y,V)+
Q_{(C)}^5(u; H,Y,V)+
Q_{(A,E)}^6(u; H,Y,V)
\}\\
 &+& \\
&& v\{ 
 \left(\n_HD_u\right)(H,Y)
+\left(\n_HF_u\right)(Y,V)\\ \nb
&& +Q_{(B,B,C, D,F)}^1(u; H,Y,V)+
Q_{(D,F,F)}^2(u; H,Y,V)+
Q_{(C,D,E)}^3(u; H,Y,V)\\ \nb
&&+Q_{(A,D,D,F)}^4(u; H,Y,V)+
Q_{(D)}^5(u; H,Y,V)+
Q_{(B,F)}^6(u; H,Y,V)
\},
\eeq
for any $Y\in T_xM$ where the horizontal lift and vertical lift are taken
at $(x,u)$.
\section{Osserman $g$-natural tangent bundles of Riemannian surfaces}\label{s4}

Let $(M,g)$ be a connected Riemannian surface, $x \in M$ and $(U, (x_1,x_2))$ a normal coordinates system on $(M,g)$ centred at $x$.
\\
For any vector $X=X^1\p_{x_1}+X^2\p_{x_2}\in T_xM$,
let us set 
\beq
{\bf i}X = -X^2\p_{x_1}+X^1\p_{x_2}.
\eeq
Then the Riemannian curvature is given by:
\beq
R(X,Y)Z = k(x)g({\bf i}X,Y){\bf i}Z
\eeq
for all vectors $X,Y,Z\in T_xM$, where $k$ denotes the Gaussian curvature of $(M,g)$. \\ \\
We have the following result:

\bpro\label{l1s4}
Let $H\in T_xM$ such that  $H_{(x,0_x)}^h$ is a unit tangent vector in \\ $(T_{(x,0_x)}TM, G_{(x,0_x)})$.
Then the spectrum of the Jacobi operator $\bar{J}_{H_{(x,0_x)}^h}$ is given by the set
\beq
\{0,\, \frac{k(x)}{(\al_1+\al_3)(0)},\ -\frac{f_6^B+k(x)(f_1^B+f_2^B)(0)}{(\al_1+\al_3)(0)}, \ 
-\frac{(f_4^B+f_5^B+f_6^B)(0)}{(\al_1+\al_3)(0)}  \}
\eeq
\epro
\noindent \ti{Proof}
\noindent Since $H\neq 0_x$,   $(H_{(x,0_x)}^h,\, 
({\bf i}H)_{(x,0_x)}^h,\, H_{(x,0_x)}^v,\, 
({\bf i}H)_{(x,0_x)}^v\, )$ is a basis in $(T_{(x,0_x)}TM$ and according to
(\ref{e1s3}) and (\ref{e2s3}) we have
\beq
\bar{J}_{H_{(x,0_x)}^h}\left( H_{(x,0_x)}^h\right)= 0_x
\eeq
\beq
\bar{J}_{H_{(x,0_x)}^h}\left(({\bf i} H)_{(x,0_x)}^h\right)= 
\frac{k(x)}{(\al_1+\al_3)(0)}({\bf i} H)_{(x,0_x)}^h
\eeq
\beq
\bar{J}_{H_{(x,0_x)}^h}\left( H_{(x,0_x)}^v\right)&=& 
-(f_4^A+f_5^A+f_6^A)(0) H_{(x,0_x)}^h\\ \nb
&&-(f_4^B+f_5^B+f_6^B)(0)H_{(x,0_x)}^v
\eeq
\beq
\bar{J}_{H_{(x,0_x)}^h}\left(({\bf i} H)_{(x,0_x)}^v\right)&=& 
-[f_6^A(0)+k(x)(f_1^A+f_2^A)(0)]({\bf i} H)_{(x,0_x)}^h\\ \nb
&&-[f_6^B(0)+k(x)(f_1^B+f_2^B)(0)]({\bf i} H)_{(x,0_x)}^v.
\eeq
Then the matrix of the operator $\bar{J}_{H_{(x,0_x)}^h}$ in the basis\\
$(H_{(x,0_x)}^h,\, 
({\bf i}H)_{(x,0_x)}^h,\, H_{(x,0_x)}^v,\, 
({\bf i}H)_{(x,0_x)}^v\, )$ is
\beq\label{e1s4}
\left(
\begin{array}{cccc}
0 & 0 & -\frac{\delta^A(0)}{(\al_1+\al_3)(0)}& 0\\
  &    &                                     & \\
0 & \frac{k(x)}{(\al_1+\al_3)(0)}&  0&-\frac{\eta^A(0)}{(\al_1+\al_3)(0)}\\
&    &                                     & \\
0& 0 & -\frac{\delta^B(0)}{(\al_1+\al_3)(0)} & 0\\
&    &                                     & \\
0 & 0 & 0 & -\frac{\eta^B(0)}{(\al_1+\al_3)(0)}
\end{array}
\right)
\eeq
where we set 
\beq
\delta^P(0)&=&(f_4^P+f_5^P+f_6^P)(0)\\
\eta^P(0)&=&-f_6^P(0)+k(x)(f_1^P+f_2^P)(0)
\eeq
for $P=A,\, B$.
This is a triangular matrix and then we get the result. $\hfill{\square}$ \\ \\
Similary arguments and Proposition 4.1 lead to the following:
\bcor
Let $\dim M = 2$. If $(TM,G)$ is a pointwise Riemannian Osserman manifold then $(M,g)$ has constant Gauss curvature.
\ecor
\noindent \ti{Proof}\\
Let $x\in M$ and $V$ a vector in $T_xM$ such that $g(V,V)=\frac{1}{\al_1(0)}$.
Then $V_{(x,0_x)}^v$ is a unit vector in $(T_{(x,0_x)}TM, G_{(x,0_x)})$ and 
  $(V_{(x,0_x)}^h,\, 
({\bf i}V)_{(x,0_x)}^h,\, V_{(x,0_x)}^v,\, 
({\bf i}V)_{(x,0_x)}^v\, )$ is a basis of $T_{(x,0_x)}TM$.\\ By computing as in the proof of Proposition \ref{l1s4} the matrix of the Jacobi operator $\bar{J}_{V_{(x,0_x)}^v}$ in this
 basis, we get:
\beq\label{e2s4}
\left(
\begin{array}{cccc}
\frac{\delta^C(0)}{\al_1(0)} & 0 &0 & 0\\
  &    &                                     & \\
0 & \frac{f_6^C(0)}{\al_1(0)}&0 & \frac{(f_6^E-f_7^E)(0)}{\al_1(0)}\\
&    &                                     & \\
\frac{\delta^D(0)}{\al_1(0)}& 0 & 0 & 0\\
&    &                                     & \\
0 & \frac{f_6^D(0)}{\al_1(0)} & 0 & \frac{(f_6^F-f_7^F)(0)}{\al_1(0)}
\end{array}
\right),
\eeq
where we put 
\beq
\delta^P(0)&=&(f_4^P+f_5^P+f_6^P)(0)\;,\; P = C , D\;.
\eeq
Hence if $(TM,G)$ is pointwise Riemannian Osserman manifold, according to Proposition \ref{l1s4}, the quotient  $\frac{k(x)}{(\al_1+\al_3)(0)}$ is necessarly an eigenvalue of the matrix (\ref{e2s4}) that does not depend on $x$. So the Gaussian curvature $k$ is constant. This completes the proof. $\hfill{\square}$
\\ \\
Let us consider the orthonormal frame bundle $\mathcal{O} (M)$ over $(M,g)$. It is a subbundle of the tangent bundle $TM$, and a $g$-natural metric on $\mathcal{O} (M)$ is the restriction of some $g$-natural metric on $TM$. It has been proved by the authors in \cite{KS2} that if $(M,g)$ has constant sectional curvature, then is orthonormal frame bundle equipped with a $g$-natural metric is  
 always locally homogeneous (cf. Corollary 4.5 in \cite{KS2}). \\
 From this observation and Proposition 4.1 we get: 
 \bcor
 Let $(M,g)$ a connected Riemannian surface, and $\tilde G$ a \\ $g$-natural metric on its orthonormal frame bundle $\mathcal{O} (M)$. Then \\$(\mathcal{O} (M) , \tilde G)$ is globally Osserman if and only if it is pointwise Osserman.
 \ecor
\noindent \ti{Proof} \\
If $(\mathcal{O} (M) , \tilde G)$ is pointwise Osserman, then by Corollary 4.1, $(M,g)$ is of cons\-tant Gaussian curvature and by Corollary 4.5 in \cite{KS2}, $(\mathcal{O} (M) , \tilde G)$ is locally homogenous. Hence the spectrum of its Jacobi operators is the same for all points and then $(\mathcal{O} (M) , \tilde G)$ is globally Osserman. $\hfill{\square}$ 
\\ \\ 
In the sequel we assume that $(M,g)$ is of constant Gaussian curvature $k$. \\ It holds:
\bpro\label{l2s4}
Let $(M,g)$ be a connected Riemannian surface with cons\-tant Gaussian curvature and $(x,u)\in TM$ with $u\neq 0_x$. Put $t = g(u,u)$. \\Then the family $(u^h,({\bf i}u)^h,u^v,({\bf i}u)^v )$ is a basis of $T_{(x,u)}TM$ and the \\non-vanishing entries of the matrix $(J_{ij})_{1\leq i,j\leq 4}$ of the Jacobi operator $\bar{J}_{u_{(x,u)}^h}$ with respect to this basis are: 
\beq
\begin{array}{lcl}
J_{22} &=& t^2\{(f_5^A-kf_1^A)[(f_4^A-kf_2^A)-(f_4^A+f_5^A+f_6^A+tf_7^A)]\\
       &&+   (f_4^C-kf_2^C)(f_5^B+k(1-f_1^B))\\
      &&-(f_5^C-kf_1^C)(f_4^B+f_5^B+f_6^B+tf_7^B)\} +kt\\
J_{42} &=& t^2\{ (f_5^A-kf_1^A)(f_4^B-kf_2^B)
  -(f_4^A+f_5^A+f_6^A+tf_7^A)(f_5^B-kf_1^B)\\
  &&+(f_4^D-kf_2^D)(f_5^B+k(1-f_1^B))
  -(f_4^B+f_5^B+f_6^B+tf_7^B)(f_5^D-kf_1^D)\}
\end{array}
\eeq
\beq
\begin{array}{lcl}
J_{13}&=& t^2[(f_4^C+f_5^C+f_6^C+tf_7^C)(f_4^D+f_5^D+f_6^D+tf_7^D)\\
     &&-(f_4^B+f_5^B+f_6^B+tf_7^B)(f_4^E+f_5^E+f_6^E+tf_7^E)]\\
   &&-t[(f_4^A+f_5^A+f_6^A+3tf_7^A)+2t({f_4^A}^\pr+{f_5^A}^\pr 
+{f_6^A}^\pr+t{f_7^A}^\pr )],\\ 
J_{33}&=& t^2\{(f_4^B+f_7^B+f_6^B+tf_7^B)
[(f_4^C+f_5^C+f_6^C+tf_7^C)-(f_4^F+f_5^F+f_6^F+tf_7^F)]\\ 
 &&+(f_4^D+f_5^D+f_6^D+tf_7^D)[(f_4^D+f_5^D+f_6^D+tf_7^D)
-(f_4^A+f_5^A+f_6^A+tf_7^A)]\}\\
&&-t[(f_4^B+f_5^B+f_6^B+3tf_7^B)+2t({f_4^B}^\pr+{f_5^B}^\pr 
+{f_6^B}^\pr+t{f_7^B}^\pr )],\\
\end{array}
\eeq
\beq
\begin{array}{lcl}
J_{24} &=& t^2\{(f_4^C-kf_2^C)[(f_4^A-kf_2^A)+(f_4^D-kf_2^D)]\\
      &&-(f_4^C-kf_2^C)(f_4^A+f_5^A+f_6^A+tf_7^A)\\
  &&-(f_5^E-kf_1^E)(f_4^B+f_5^B+f_6^B+tf_7^B)\}\\
 && -t[(f_6^A+tf_7^A)+k(f_1^A+f_2^A)],\\
J_{44} &=& t^2\{ (f_4^D-kf_2^D)^2+(f_4^B-kf_2^B)(f_4^C-kf_2^C)\\
       &&-(f_4^D-kf_2^D)(f_4^A+f_5^A+f_6^A+tf_7^A)
-(f_5^F-kf_1^F)(f_4^B+f_5^B+f_6^B+tf_7^B)\}\\
&&-t[(f_6^B+tf_7^B)+k(f_1^B+f_2^B)]\;.\\
\end{array}
\eeq

\epro
\brmk\label{r1s4}
\ben
\item It is easy to check that
\beq
(\phi_1+\phi_3)J_{13}+\phi_2 J_{33} = 0,
\eeq
\beq
\al_2(J_{44}-J_{22})+(\al_1+\al_3)J_{24} = \al_1J_{42}.
\eeq
\item The following vectors 
\beq
v_1 &=& \frac{1}{\sqrt{t(\phi_1+\phi_3)(t)}}\ u^h,\\ \nb
            & & \\ 
v_2 &=& \sqrt{\frac{(\phi_1+\phi_3)(t)}{t\phi(t)}}\ u^v -\frac{\phi_2(t)}{\sqrt{t\phi(t)(\phi_1+\phi_3)(t)}}\ u^h,\\ \nb
         & & \\ 
v_3 &=& \frac{1}{\sqrt{t(\al_1+\al_3)(t)}}\ ({\bf i}u)^h,\\ \nb
       & & \\ 
v_4 &=& \sqrt{\frac{(\al_1+\al_3)(t)}{t\al(t)}}\ ({\bf i}u)^v -\frac{\al_2(t)}{\sqrt{t\al(t)(\al_1+\al_3)(t)}}\ ({\bf i}u)^h,
\eeq
where the lifts are taken at (x,u),
determine an orthonormal basis of \\
$(T_{(x,u)}TM,\, G_{(x,u)})$.
\een
\ermk
\bpro\label{s4l3}
Let $(x,u)\in TM$ such that $u\neq 0_x$ and $t = g(u,u)$.  \\ Then the spectrum  of  Jacobi operator $\bar{J}_{u_{(x,u)}^h}$ is given by the set
\beq
\{0,\, J_{33},\, \frac{(J_{22}+J_{44})+\sqrt{\Delta}}{2},\, 
\frac{(J_{22}+J_{44})-\sqrt{\Delta}}{2}\},
\eeq
where $\Delta= \left(J_{22}-J_{44}+2\frac{\al_2}{\al_1+\al_3}J_{42}\right)^2
+4\frac{\al}{(\al_1+\al_3)^2}J_{42}^2$.
\epro
\ti{Proof}\\
According to Remark \ref{r1s4} and Proposition \ref{l2s4} the matrix of $\bar{J}_{u_{(x,u)}^h}$ in the orthonormal basis 
$(v_1,\, v_2,\, v_3,\, v_4)$ is given by 
\beq
\left(
\begin{array}{cccc}
0 & 0 &0 & 0\\
  &    &                                     & \\
0 & J_{33}&0 & 0\\
&    &                                     & \\
0& 0 & (J_{22}+\frac{\al_2}{\al_1+\al_3}J_{42}) & \frac{\sqrt{\al}}{\al_1+\al_3}J_{42}\\
&    &                                     & \\
0 & 0 &  \frac{\sqrt{\al}}{\al_1+\al_3}J_{42}& (J_{44}-\frac{\al_2}{\al_1+\al_3}J_{42})
\end{array}
\right).
\eeq
So by computing the eigenvalues of this matrix, we obtain the proof.
$\hfill{\square}$
\\ \\
Using Proposition \ref{s4l3} and by notifying that $G(u^h,u^h)=t(\phi_1+\phi_3)(t)$ with $t=g(u,u)$, we obtain the following results:

\bthm\label{p1s4}
$(TM,G)$ is a pointwise  Osserman manifold if and only if
\ben
\item $(M,g)$ has constant Gauss curvature $k$.
\item The  eigenvalues of its 
Jacobi operators on the unit tangent  bundle $S(TTM)$   are  the functions $(\lambda_i)_{i=1,2, 3}$
 defined on $TM$ by
\beq\label{f1}
\lambda_0(x,u)&=&0,\\ \nb
\lambda_1(x,u)&=& \frac{J_{33}}{t(\phi_1+\phi_3)},\\ \nb
\lambda_2(x,u)&=& \frac{(J_{22}+J_{44})+\sqrt{\Delta}}{2t(\phi_1+\phi_3)},\\ \nb
\lambda_3(x,u) &=& \frac{(J_{22}+J_{44})-\sqrt{\Delta}}{2t(\phi_1+\phi_3)},
\eeq
if $u\neq 0_x$ \\
and
\beq\label{f2}
\lambda_0(x,0_x)&=&0,\\ \nb
\lambda_1(x,0_x)&=& -\frac{(f_4^B+f_5^B+f_6^B)(0)}{(\al_1+\al_3)(0)},\\ \nb
\lambda_2(x,0_x)&=& \frac{k}{(\al_1+\al_3)(0)},\\ \nb
\lambda_3(x,0_x) &=&-\frac{f_6^B+k(f_1^B+f_2^B)(0)}{(\al_1+\al_3)(0)} .
\eeq. 
\een
\ethm
\bthm
  $(TM,G)$ is a globally  Osserman manifold if and only if 
\ben
\item $(M,g)$ has constant Gauss curvature $k$.
\item  The  eigenvalues of its Jacobi operators on the unit tangent bundle $S(TTM)$   are  the real numbers 
$(\tilde{\lambda}_i)_{i=1,2, 3}$
 given by
\beq\label{f3}
\tilde{\lambda}_0&=&0,\\ \nb
\tilde{\lambda}_1 &=& -\frac{(f_4^B+f_5^B+f_6^B)(0)}{(\al_1+\al_3)(0)},\\ \nb
\tilde{\lambda}_2 &=& \frac{k}{(\al_1+\al_3)(0)},\\ \nb
\tilde{\lambda}_3  &=&-\frac{f_6^B+k(f_1^B+f_2^B)(0)}{(\al_1+\al_3)(0)} .
\eeq. 
\een
\ethm
\noindent In the following we apply the result in Theorem 4.2 to the Sasaki metric and to the Cheeger-Gromoll metric on the tangent bundle. \\ \\
{\bf Applications}: \\
1. Let $G$ be the Sasaki metric on the tangent bundle $TM$. \\
In this case the functions $\al_i$ and $\be_i$ of Proposition 2.1 are given by:
$$ \al_1 = 1\;;\;\al_2 = \al_3 = 0 \;\;\mbox{and}$$
$$\be_1 = \be_2 = \be_3 = 0 \;.$$ 
The eigenvalues $\tilde \lambda_0$, $\tilde \lambda_1$, $\tilde \lambda_2$, $\tilde \lambda_3$ of Theorem 4.2 are:
$$\tilde \lambda_0 = \tilde \lambda_1 = 0\;;\;\tilde \lambda_2 = k\;;\;\tilde \lambda_3 = 0\;.$$ 
\\
2. Let $G$ be the Cheeger-Gromoll metric on the tangent bundle $TM$. \\
Then the functions $\al_i$ and $\be_i$ of Proposition 2.1 are given by:
$$ \al_1 = \be_1 = \frac{1}{1+2t}\;;\;\al_2 = \al_2 = 0 \;\;\mbox{and}$$
$$\al_3 = \frac{2t}{1+2t}\;;\; \be_3 = - \frac{1}{1+2t}\;.$$ 
The eigenvalues $\tilde \lambda_0$, $\tilde \lambda_1$, $\tilde \lambda_2$, $\tilde \lambda_3$ of Theorem 4.2 are in this case:
$$\tilde \lambda_0 = \tilde \lambda_1 = 0\;;\;\tilde \lambda_2 = k\;;\;\tilde \lambda_3 = 0\;.$$
\\
We can conclude that the tangent bundle $TM$ with the Sasaki metric or the Cheeger-Gromoll metric is globally Ossermann if and only if $(M,g)$ is of cons\-tant Gaussian curvature $k$ and the eigenvalues of its Jacobi operators are $0$ (with multiplicity three) and $k$. 
\\ \\
\noindent The following consequence for the sectional curvature of $g$-natural metrics can be derived from Theorem 4.2:
\bcor
Only flat $g$-natural metrics on the tangent bundle of a Riemannian surface $(M,g)$ are of constant sectional curvature. 
\ecor
\noindent \ti{Proof}\\
Let $G$ be a $g$-natural metric on $TM$ of constant sectional curvature. Then $(TM , G)$ is globally Osserman. But also $(M,g)$ is  then flat (cf. \cite{DET}) and according to (\ref{f3}), the eigenvalue $\bar{\lambda}_2 = \frac{k}{(\al_1+\al_3)(0)}$ of the Jacobi operators of $(TM,G)$ is like the eigenvalue $\bar{\lambda}_0$ equal to zero. \\Thus $0$ is an eigenvalue of the Jacobi operators of $(TM,G)$ with multiplicity  at least two. Hence $(TM,G)$ is flat. $\hfill{\square}$

\brmk
This corollary extends  \ti{Proposition 4.3} in \cite{DET} to the case where $\dim M =2$.
\ermk

\bibliographystyle{plain}

\end{document}